\documentstyle[12pt]{article}

\newtheorem{lemma}{Lemma}
[section]

\newtheorem{thm}{Theorem}
[section]
\newtheorem{cor}{Corollary}
[section]
\newtheorem{rmk}{Remark}
[section]

\newcommand{\mysection}[1]{\section{#1}\setcounter{equation}{0}}

\newcommand{\la}{\lambda}

\newcommand{\ptl}{\partial}

\newcommand{\ra}{\rightarrow}
\newcommand{\lra}{\longrightarrow}

\def\a{\alpha}

\def\R2{{\Bbb R}^2}
\def\Bbb{\mathbb}

\def\beq{\begin{equation}}
\def\eeq{\end{equation}}
\def\ba{\begin{array}}
\def\ea{\end{array}}
\def\barr{\begin{eqnarray}}
\def\earr{\end{eqnarray}}

\def\part{\partial}
\def\fr{\frac}

\newcommand{\ls}{\setlength{\baselineskip}{18pt}
                      \setlength{\parskip}{3mm} }

\title {The $C^{\a}$ regularity of a class of non-homogeneous ultraparabolic equations}
\author{ WANG Wendong and ZHANG Liqun
\thanks{The author currently is working at NSFC.The research is partially supported by the
Chinese NSF under grant 10325104. Email: wangwendong@amss.ac.cn and
lqzhang@math.ac.cn} }

\date{Institute of Mathematics, AMSS, Academia Sinica,
Beijing\\  }

\begin{document}

\maketitle

\begin{abstract}
We obtain the $C^{\a}$ regularity for weak solutions of a class of
non-homogeneous ultraparabolic equation, with measurable
coefficients. The result generalizes our recent $C^{\a}$ regularity
results of homogeneous ultraparabolic equations.

\end{abstract}

{\small keywords: Non-homogeneous, ultraparabolic equations,
$C^{\a}$ regularity }

\pagenumbering{arabic}
%%%%%%%%%%%%%%%%%%%%%%%%%%%%%%%%%%%%%%%%%%%%%%%%%%%%%%%%%%%%%%%%%%%%%%%%%%%%%%%%
\mysection{Introduction} \ls \noindent

The regularity of ultraparabolic equation becomes important since it
has many applications.  From mathematical points of view, it has
some special algebraic structures and is degenerated. Though there
are more and more studies on this problem in recent years, it is
still unclear in general, whether the interior $C^{\a}$ regularity
results hold for weak solutions of the ultraparabolic equations with
bounded measurable coefficients like the parabolic cases.

One of the typical example of the ultraparabolic equation is the
following equation
$$
\displaystyle \fr{\partial \,u}{\partial \,t}+y \fr{\partial
\,u}{\partial \,x}-u^2\fr{\partial^2 \,u}{\partial
\,y^2}=0.\leqno(1.1)
$$
This is of strong degenerated parabolic type equations, more
precisely, an ultraparabolic type equation. However, if the
coefficient is smooth it satisfies the well known H\"ormander's
hypoellipticity conditions, which sheds lights on the smoothness of
weak solutions. It is interesting if the weak solution of equation
(1.1) is still smooth when the coefficient is only measurable
functions.

On the other hand, the equation (1.1), if consider it as an equation
of $\fr 1 u$, has the divergent form. A recent paper by Pascucci and
Polidoro [12], Cinti, Pascucci and Polidoro [2] proved that the
Moser iterative method still works for a class of ultraparabolic
equations with measurable coefficients. Their results show that for
a non-negative sub-solution $u$ of (1.1), the $L^{\infty}$ norm of
$u$ is bounded by the $L^p$ norm ($p \ge 1$). This is a very
important step to the final regularity of solutions of the
ultraparabolic equations.

We seems to have proved in [15], [17] that the weak solution
obtained in [14] of (1.1) is of $C^{\a}$ class, then $u$ is smooth.
In this paper, we are concerned with the $C^{\a}$ regularity of
solutions of more general ultraparabolic equations.

We consider a class of non-homogeneous Kolmogorov-Fokker-Planck type
operator on ${R}^{N+1}$:
$$
\displaystyle  Lu \equiv
\sum_{i,j=1}^{m_0}{\ptl_{x_i}(a_{ij}(x,t)\ptl_{x_j}\,u
)}+\sum_{i,j=1}^N b_{ij}x_i {\ptl_{x_j}\, u }- {\ptl_t
\,u}=0,\leqno(1.2)
$$
where $(x,t)\in { R}^{N+1}$, $1\leq m_0\leq N$, and $b_{ij}$ is
constant for every $i,j= 1,\cdots, N$. Let $A=(a_{ij})_{N\times N},$
where $a_{ij}=0,$ if $i>m_0$ or $j>m_0.$ We make the following
assumptions on the coefficients of $L$:

$(H_1)$  $a_{ij}=a_{ji} \in L^{\infty} ({R}^{N+1})$ and there
exists a $\lambda >0$ such that
$$
\fr{1}{\lambda}\sum_{i=1}^{m_0}\xi_i^2 \leq \sum_{i,j=1}^{m_0}
a_{ij}(x,t)\xi_i \xi_j \leq {\lambda}\sum_{i=1}^{m_0}\xi_i^2
$$
for every $(x,t)\in {R}^{N+1}$, and $\xi \in {R}^{m_0}$.

$(H_2)$ The matrix $B=(b_{ij})_{N \times N}$ has the form
$$\left(
\begin{array}{ccccc}
\ast & {B_1} & 0 & \cdots & 0  \\
\ast &  \ast  & {B_2} & \cdots & 0 \\
\vdots & \vdots & \vdots & \ddots &\vdots \\
\ast &  \ast & \ast & \cdots & {B_d} \\
\ast &  \ast & \ast & \cdots & \ast
\end{array}
\right)
$$
where $B_k$ is a matrix $m_{k-1}\times m_{k}$ with rank $m_k$ and
$m_0\geq m_1\geq \cdots\geq m_d$, $m_0+m_1+\cdots+m_d=N$.

The requirements of matrix $B$ in $(H_2)$ ensure that the operator
$L$ with the constant $a_{ij}$ satisfies the well-known
H\"ormander's hypoellipticity condition. We let $\la$ satisfies
$||B||\leq \la$ where the norm $||\cdot||$ is in the sense of matrix
norm. We refer [2] for more details on non-homogeneous
Kolmogorov-Fokker-Planck type operator on ${R}^{N+1}$.

The Schauder type estimate of (1.2) has been obtained for example,
in [18], [19] and [16]. Besides, the regularity of weak solutions
have been studied by Bramanti, Cerutti and Manfredini [1], Polidoro
and Ragusa [13] assuming a weak continuity on the coefficient
$a_{ij}$. It is quite interesting whether the weak solution has
H\"older regularity under the assumption $(H_1)$ on $a_{ij}$. One of
the approach to the H\"older estimates is to obtain the Harnack type
inequality. In the case of elliptic equations with measurable
coefficients, the Harnack inequality is obtained by J. Moser [9] via
an estimate of BMO functions due to F. John and L. Nirenberg
together with the Moser iteration method. J. Moser [10] also
obtained the Harnack inequality for parabolic equations with
measurable coefficients by generalizing the John-Nirenberg estimates
to the parabolic case. Another approach to the H\"older estimates is
given by S. N. Kruzhkov [6], [8] based on the Moser iteration to
obtain a local priori estimates, which provides a short proof for
the parabolic equations. Nash [11] introduced another technique
relying on the ${\rm Poincar\acute{e}}$ inequality and obtained the
H\"older regularity. Also De Giorgi developed an approach to obtain
the H\"older regularity for elliptic equations.

We prove a ${\rm Poincar\acute{e}}$ type inequality for non-negative
weak sub-solutions of (1.2). Then we apply it to obtain a local
priori estimates which implies the H\"older estimates for
ultraparabolic equation (1.2).

Let $D_{m_0}$ be the gradient with respect to the variables $x_1,
x_2,\cdots, x_{m_0}$. And
$$Y=\sum_{i,j=1}^N b_{ij}x_i {\ptl_{x_j}}- {\ptl_t}.$$
We say that $u$ is a $weak \,solution$ if it satisfies (1.2) in the
distribution sense, that is for any $\phi \in C^{\infty}_0(\Omega)$,
where $\Omega$ is a open subset of $R^{N+1}$, then
$$
\int_{\Omega} \phi Yu-(Du)^T AD\phi = 0, \leqno(1.3)
$$
and $u$, $D_{m_0}u$, $Yu \in L^2_{\rm loc}(\Omega).$

Our main result is the following theorem.
\begin{thm}
Under the assumptions $(H_1)$ and $(H_2)$, the weak solution of
(1.2) is H\"older continuous.
\end{thm}

\mysection{Some Preliminary Results}

One of the important feature of equation (1.2) is that the
fundamental solution can be written explicitly if the coefficients
$a_{ij}$ is constant (cf. [4], [7]). Besides, there are some
geometric and algebraic structures in the space $R^{N+1}$ induced by
the constant matrix $B$ (see for instance, [7]).

We follow the earlier notations and give some basic properties used
for example, by [2] and [7], and more details see [2] and [7].

Let $E(\tau)=\rm{exp}(-\tau B^T)$. For $(x,t), (\xi,\tau) \in
R^{N+1}$, set
$$(x,t)\circ (\xi,\tau)=(\xi+E(\tau)x,t+\tau),$$ then $(R^{N+1}, \circ)$
is a Lie group with identity element $(0,0)$, and the inverse of an
element is $(x,t)^{-1}=(-E(-t)x,-t)$. The left translation by
$(\xi,\tau)$ given by
$$(x,t)\mapsto (\xi,\tau)\circ (x,t),$$
is a invariant translation to operator $L$ when coefficient $a_{ij}$
is constant.
The associated dilation to operator $L$ with constant
coefficient $a_{ij}$ is given by
$$
\delta_t=diag (t I_{m_0},t^3 I_{m_1},\cdots,t^{2d+1}I_{m_d},t^2),
$$
where $I_{m_k}$ denotes the $m_k\times m_k$ identity matrix, $t$ is
a positive parameter, also we assume
$$
D_t=diag(t I_{m_0},t^3 I_{m_1},\cdots,t^{2d+1}I_{m_d}),
$$
and denote
$$Q=m_0+3m_1+\cdots+(2d+1)m_d,$$
then the number $Q+2$ is usually called the homogeneous dimension of
$(R^{N+1},\circ)$ with respect to the dilation $\delta_t$.

The norm in $R^{N+1}$, related to the group of translations and
dilation to the equation is defined by $$||(x,t)||=r,$$ if $r$ is
the unique positive solution to the equation
$$
\fr{x_1^2}{r^{2\a_1}}+\fr{x_2^2}{r^{2\a_2}}+\cdots+\fr{x_N^2}{r^{2\a_N}}
+\fr{t^2}{r^4}=1,
$$
where $(x,t) \in R^{N+1}\setminus \{0\}$ and
$$
\a_1=\cdots=\a_{m_0}=1, \quad
\a_{m_0+1}=\cdots=\a_{m_0+m_1}=3,\cdots,
$$
$$
\a_{m_0+\cdots+m_{d-1}+1}=\cdots=\a_N=2d+1.
$$
And $||(0,0)||=0$. The balls at a point $(x_0,t_0)$ is defined by
$${\cal B}_r(x_0,t_0)=\{(x,t)|\quad ||(x_0,t_0)^{-1}\circ (x,t)||\leq r\},$$
and
$${\cal B}^-_r(x_0,t_0)={\cal B}_r(x_0,t_0)\cap\{t<t_0\}.$$
For convenience, we sometimes use the cube replace the balls. The
cube at point $(0,0)$ is given by
$$
{\cal C}_r(0,0)=\{(x,t)|\quad |t|\leq r^2,\quad |x_1|\leq r^{\a_1},
\cdots, |x_N|\leq r^{\a_N}\}.
$$
It is easy to see that there exists a constant $\Lambda$ such that
$$
{\cal C}_{\fr r \Lambda}(0,0)\subset{\cal B}_r(0,0)\subset{\cal
C}_{\Lambda r}(0,0),
$$
where $\Lambda$ only depends on $B$ and $N$.

When the matrix $(a_{ij})_{N \times N}$ is of constant matrix, we
denoted it by $A_0$, and $A_0$ has the form
$$
A_0= \left(
\begin{array}{cc}
I_{m_0}& 0  \\
0 & 0
\end{array}
\right)
$$
then let
$$
{\mathcal C}(t)\equiv\int_0^t E(s)A_0 E^T(s)ds,
$$
which is positive when $t>0,$ and  the operator $L_1$ takes the form
$$
L_1=div (A_0D)+Y,
$$
whose fundamental solution $\Gamma_1(\cdot,\zeta)$ with pole in
$\zeta\in R^{N+1}$ has been constructed as follows:
$$
\Gamma_1(z,\zeta)=\Gamma_1(\zeta^{-1}\circ z,0), \qquad z, \zeta \in
R^{N+1},\quad z \neq \zeta,
$$
where $z=(x,t)$. And $\Gamma_1(z,0)$ can be written down explicitly
$$
\Gamma_1(z,0)= \Big {\{}
\begin{array}{cc}{\frac{(4\,\pi)^{-\frac{N}{2}}}{\sqrt{\det {\mathcal C}(t)}}\exp(-\frac{1}{4}\langle {\mathcal C}
^{-1}(t)x,x\rangle -t \, tr(B))} & {\rm if} \quad t>0,
\\
0 & {\rm if} \quad t\leq 0. \end{array}\leqno(2.1)
$$
There are some basic estimates for $\Gamma_1$ (see [2])
$$
\Gamma_1(z,\zeta)\leq C ||\zeta^{-1}\circ z||^{-Q}, \leqno(2.2)
$$
$$
|\ptl_{\xi_i}\,\Gamma_1(z,\zeta)|\leq C ||\zeta^{-1}\circ
z||^{-Q-1}, \leqno(2.3)
$$
where $i=1,\cdots,m_0$, for all $z,\zeta\in R^N\times(0,T]$.

A weak sub-solution of (1.2) in a domain $\Omega$ is a function
$u$ such that $u$, $D_{m_0}u$, $Yu \in L^2_{loc}(\Omega)$ and for
any $\phi \in C^{\infty}_0(\Omega)$, $\phi \geq 0$,
$$
\int_{\Omega} \phi Yu-(Du)^T AD\phi \geq 0. \leqno(2.4)
$$
Similarly, let $Y_0=<x,\,B_0 D>-\partial_ t$, where $B_0 \,$ has the
form
$$\left(
\begin{array}{cccccc}
0 & {B_1} & 0 & \cdots & 0 \\
0 & 0  & {B_2} & \cdots & 0 \\
\vdots & \vdots & \vdots & \ddots &\vdots \\
0 & 0 & 0 & \cdots & {B_d} \\
0 & 0 & 0 & \cdots & 0
\end{array}
\right)
$$
We denote $L_0=div (A_0 D)+Y_0,$ and can define in the same way
$E_0(t)$, ${\mathcal C}_0(t),$ and $\Gamma_0(z,\zeta)$ with respect
to $B_0.$ We recall that ${\mathcal C}_0(t) (t>0)$ (see[7])
satisfies
$$
{\mathcal C}_0(t)=D_{t^{\fr 1 2}}{\mathcal C}_0(1)D_{t^{\fr 1 2}}.
\leqno(2.5)
$$

The following lemma is obtained by Lanconelli and Polidoro (see
[7]), which is need in our proof.
\begin{lemma} In addition to the above assumptions,  for every given
$T>0$, there exist positive constants $C_T$ and $C'_T$ such that
$$
\langle {\mathcal C}_0(t)x,x\rangle(1-C_T\,t)\leq \langle {\mathcal
C}(t)x,x\rangle \leq \langle {\mathcal
C}_0(t)x,x\rangle(1+C_T\,t),\leqno(2.6)
$$
$$
\langle {\mathcal C}_0^{-1}(t)x,x\rangle(1-C_T\,t)\leq \langle
{\mathcal C}^{-1}(t)x,x\rangle \leq \langle {\mathcal
C}_0^{-1}(t)x,x\rangle (1+C_T\,t),\leqno(2.7)
$$
$$
C_T^{'-1}t^Q(1-C_T\,t)\leq {\rm det} {\mathcal C}(t)\leq
C'_T\,t^Q(1+C_T\,t),\,\leqno(2.8)
$$
for every $(x,t)\in R^N\times (0,T]\,$ and
$t<\,\fr{1}{C_T}.$ \\

\end{lemma}

A result of Cinti, Pascucci and Polidoro obtained by using the
Moser's iterative method (see [2]) states as follows.

\begin{lemma}
Let $u$ be a non-negative weak sub-solution of (1.2) in $\Omega$.
Let $(x_0,t_0)\in \Omega$ and $\overline{{\cal
B}^-_r(x_0,t_0)}\subset \Omega$ and $p \geq 1$. Then there exists a
positive constant $C$ which depends only on $\la$ and $Q$ such that,
for $0 < r\leq 1$
$$
\sup_{{\cal B}^-_{\fr r 2}(x_0,t_0)} u^p \leq \fr
{C}{r^{Q+2}}\int_{{\cal B}^-_r(x_0,t_0)} u^p,\leqno(2.9)
$$
provided that the last integral converges.
\end{lemma}

We copy a classical potential estimates (cf. (1.11) in [3]) here to
prove the ${\rm Poincar\acute{e}}$ type inequality.

\begin{lemma}
Let $(R^{N+1},\circ)$ is a homogeneous Lie group of homogeneous
dimension $Q+2$, $\a \in (0, Q+2)$ and $G \in C(R^{N+1}\setminus
\{0\})$ be a $\delta_{\la}$-homogeneous function of degree $\a-Q-2$.
If $f \in L^p(R^{N+1})$ for some $p \in (1,\infty)$, then
$$
G_f(z)\equiv \int_{R^{N+1}} G(\zeta ^{-1}\circ z)f(\zeta)d\zeta,
$$
is defined almost everywhere and there exists a constant
$C=C(Q,p)$ such that
$$
||G_f||_{L^q(R^{N+1})}\leq C \max_{||z||=1} |G(z)|\quad
||f||_{L^p(R^{N+1})},\leqno(2.10)
$$
where $q$ is defined by
$$
\fr 1q =\fr 1p-\fr{\a}{Q+2}.
$$
\end{lemma}
\begin{cor} Let $f\in L^2(R^{N+1})$,  recall the definitions in [2]
$$
\Gamma_1(f)(z)=\int_{R^{N+1}}\Gamma_1(z,\zeta)f(\zeta) d\zeta,
\qquad \forall z\in R^{N+1},
$$
and
$$
\Gamma_1(D_{m_0}f)(z)=-\int_{R^{N+1}}D_{m_0}^{(\zeta)}\Gamma_1(z,\zeta)f(\zeta)
d\zeta, \qquad \forall z\in R^{N+1},
$$
then exists a positive constant $C=C(Q,T,B)$ such that
$$
\|\Gamma_1(f)\|_{L^{2\tilde{k}}(S_T)}\leq C\|f\|_{L^2(S_T)},
\leqno(2.11)
$$
and
$$
\|\Gamma_1(D_{m_0}f)\|_{L^{2k}(S_T)}\leq C\|f\|_{L^2(S_T)},
\leqno(2.12)
$$
where $\tilde{k}=1+\fr{4}{Q-2}$, $k=1+\fr{2}{Q}$ and
$S_T=R^N\times]0,T].$
\end{cor}

\mysection{Proof of Main Theorem}

To obtain a local estimates of solutions of the equation (1.2), for
instance, at point $(x_0,t_0)$, we may consider the estimates at a
ball centered at $(0,0)$, since the equation (1.2) is invariant
under the left group translation when $a_{ij}$ is constant. By
introducing a ${\rm Poincar\acute{e}}$ type inequality, we prove the
following Lemma 3.5 which is essential in the oscillation estimates
in Kruzhkov's approaches in parabolic case. Then the $C^{\a}$
regularity result follows easily by the standard arguments.

For convenience, in the following discussion, we let
$x'=(x_1,\cdots,x_{m_0})$ and $x=(x', \overline x)$. We consider the
estimates in the following cube, instead of ${\cal B}^-_r$,
$$
{\cal C}_r^{-}=\{(x,t)| \quad-r^2\leq t \leq 0, |x'|\leq r,
|x_{m_0+1}|\leq (\la N^2 r)^{3}, \cdots, |x_N|\leq (\la N^2
r)^{2d+1}\}.
$$
Let
$$
K_r=\{x'|\quad |x'|\leq r \},
$$
$$
S_r=\{ \overline x\quad|\quad|x_{m_0+1}|\leq (\la N^2 r)^{3},
\cdots, |x_N|\leq (\la N^2 r)^{2d+1}\}.
$$

Let $0<\a, \beta<1$ be constants, for fixed $t$ and $h$, let
$$
{\cal N}_{t,h}=\{(x',\overline x)\in K_{\beta r}\times S_{\beta
r},\quad u(\cdot,t) \geq h\}.
$$
In the following discussions, we sometimes abuse the notations of
${\cal B}^-_r$ and ${\cal C}_r^-$, since there are equivalent, and
we always assume $r \ll 1$ and $\la>8$ in the following arguments,
since $\la$ can choose a large constant. Moreover, all constants
depend on $m_0$, $d$ or $Q$ will be denoted by dependence on $B$.

\begin{lemma}
Suppose that $u(x,t)\geq 0$ be a solution of equation (1.2) in
${\cal B}^-_r$ centered at $(0,0)$ and
$$
mes\{(x,t)\in {\cal B}^-_r, \quad u \geq 1\} \geq \fr 1 2 mes ({\cal
B}^-_r).
$$
Then there exist constants $\a$, $\beta$ and $h$, $0<\a, \beta, h<1$
which only depend on $B$, $\la$ and $N$ such that for almost all
$t\in (-\a r^2,0)$,
$$
mes\{{\cal N}_{t,h}\} \geq \fr {1}{11}mes\{ K_{\beta r}\times
S_{\beta r}\}.
$$
\end{lemma}
{\it Proof:} Let
$$
v=\ln^+(\fr{1}{u+h^{\fr 9 8}}),
$$
where $h$ is a constant, $0<h<1$, to be determined later. Then $v$
at points where $v$ is positive, satisfies
$$
\displaystyle
\sum_{i,j=1}^{m_0}{\ptl_{x_i}\,(a_{ij}(x,t)\ptl_{x_j}\,v )}-(Dv)^TA
Dv+x^T B Dv - {\ptl_t \,v}=0.\leqno(3.1)
$$
Let $\eta(x')$ be a smooth cut-off function so that
$$
\eta(x')=1,\quad \hbox {for} \quad |x'|< \beta r,
$$
$$
\eta(x')=0,\quad \hbox {for} \quad |x'|\geq r.
$$
Moreover, $0\leq\eta \leq 1$ and $|D_{m_0} \eta|\leq \fr
{2m_0}{(1-\beta)r}$.

Multiplying $\eta^2(x')$ to (3.1) and integrating by parts on
$K_r\times S_{\beta r}\times(\tau,t)$
$$
\ba{lllll} \int_{K_{\beta r}}\int_{S_{\beta r}} v(t,x',\overline x)d
\overline x dx' +\fr {1}{2\la}\int_\tau^t \int_{K_{
r}}\int_{S_{\beta r}}\eta^2 \,
|D_{m_0}v|^2d \overline x dx'dt \\ \\
\leq \fr {C} {\beta^{Q}(1-\beta)^2}mes(S_{\beta r})mes(K_{\beta
r})+\int_\tau^t
\int_{K_{r}}\int_{S_{\beta r}} \eta^2 x^TBDv d \overline xdx' dt \\ \\
\quad +\int_{K_{r}}\int_{S_{\beta r}} v(\tau,x',\overline x)d
\overline x dx',\qquad a.e. \quad\tau, t\in(-r^2,0),\ea\leqno(3.2)
$$
where $C$ only depends on $\la$, $B$ and $N$. Let
$$
I_B\equiv\int_{K_{r}}\int_{S_{\beta r}} \eta^2 \sum_{i,j=1}^N
x_ib_{ij}\ptl_{x_j}v d \overline x dx'= I_{B_1}+I_{B_2},
$$
where
$$
I_{B_1}=\int_{K_{r}}\int_{S_{\beta r}} \eta^2 \sum_{i=1}^N
\sum_{j=1}^{m_0} x_ib_{ij}\ptl_{x_j}v d \overline x dx',
$$
$$
I_{B_2}=\int_{K_{r}}\int_{S_{\beta r}} \eta^2 \sum_{i=1}^N
\sum_{j=m_0+1}^{N} x_ib_{ij}\ptl_{x_j}v d \overline x dx'.
$$
On the other hand
$$
\ba{llll} |I_{B_1}|\leq\int_{K_{r}}\int_{S_{\beta r}}
\varepsilon\eta^2|D_{m_0}v|^2+C_\varepsilon\eta^2\sum_{j=1}^{m_0}\sum_{i=1}^N
|x_i b_{ij}|^2 d \overline x dx' \\ \\
\leq\int_{K_{r}}\int_{S_{\beta r}}
\varepsilon\eta^2|D_{m_0}v|^2d\overline x
dx'+C(\varepsilon,B,\lambda,N)\beta^{-Q}|K_{\beta r}||S_{\beta
r}|,\ea\leqno(3.3)
$$
and
$$
\ba{llllllll} |I_{B_2}| &\leq & |\int_{K_{r}}\int_{S_{\beta r}}
\eta^2 \sum_{i=1}^N \sum_{j=m_0+1}^{N} x_i b_{ij}\ptl_{x_j}v
d\overline x dx'| \\ \\& \leq & |\int_{K_{r}}\int_{S_{\beta r}}
{-}\eta^2\sum_{i=1}^N\sum_{j>m_0}\delta_{ij}b_{ij}vd \overline x dx'|\\
\\&&+|\int_{K_{r}}\int_{\partial_j S_{\beta r}}
\eta^2\sum_{i=1}^N\sum_{j>m_0}x_{i}b_{ij}v d\overline {x_j} dx'| \\
\\& \leq & \lambda N\beta^{-Q}|K_{\beta r}||S_{\beta r}| \ln
(h^{-\fr 9 8})\\ \\&& +\la\sum_{i=1}^N\sum_{j>m_0}\fr{(\lambda N^2
r)^{\alpha_i}}{(\lambda N^2 r)^{\alpha_j}} \beta^{-2Q}|K_{\beta
r}||S_{\beta r}| \ln (h^{-\fr 9 8}), \ea
$$
where $\overline {x_j}=(x_{m_0+1}, \dots, x_{j-1}, x_{j+1}, \dots,
x_N)$. When $\a_i\geq\a_j$, we have
$$
\int_\tau^t |I_{B_2}| \leq (\lambda N \,r^2+\la r^2
N^2)\beta^{-2Q}|K_{\beta r}||S_{\beta r}| \ln (h^{-\fr 9 8}),
$$
or $i<j$, thus $\a_j=\a_i+2$ by the property of $B$, then
$$
\int_\tau^t |I_{B_2}| \leq (\lambda N \,r^2+\la^{-1}
N^{-2})\beta^{-2Q}|K_{\beta r}||S_{\beta r}| \ln (h^{-\fr 9 8}).
$$
By $\la>8$ choose $r$ small enough, such that
$$\lambda N \,r^2+\la r^2N^2+\lambda^{-1}N^{-2}<\fr{1}{8},$$ thus
$$
\int_\tau^t |I_{B_2}|\leq \fr{1}{4}\beta^{-2Q}|K_{\beta r}||S_{\beta
r}| \ln (h^{-\fr 9 8}). \leqno (3.4)
$$
Integrating by t to $I_B$, we have
$$
\ba{llllll} \hspace*{-10pt}\int_\tau^t\int_{K_{r}}\int_{S_{\beta r}}
\eta^2 x^T B Dv d \overline x dx'dt \\ \\ \leq \fr {1}{4}
{\beta}^{-2Q}
\ln(h^{-\fr 9 8}) mes(S_{\beta r})mes(K_{\beta r})\\
\\\hspace*{-3pt}+\int_\tau^t\int_{K_{r}}\int_{S_{\beta r}}
\varepsilon\eta^2|D_{m_0}v|^2+C(\varepsilon,B,\lambda,N)\beta^{-Q}|K_{\beta
r}||S_{\beta r}|. \ea\leqno (3.5)
$$

We shall estimate the measure of the set ${\cal N}_{t,h}$. Let
$$
\mu(t)=mes\{(x',\overline x)|\quad x'\in K_r,\, \overline x \in
S_{r}, \, u(\cdot,t)\geq 1\}.
$$
By our assumption, for $0<\a< \fr 12$
$$
\fr 12 r^2 mes(S_{r})mes(K_{r})\leq \int_{-r^2}^0
\mu(t)dt=\int_{-r^2}^{-\a r^2}\mu(t)dt+\int_{-\a r^2}^{0}\mu(t)dt,
$$
that is
$$
\int_{-r^2}^{-\a r^2}\mu(t)dt\geq (\fr 12-\a)r^2
mes(S_{r})mes(K_{r}),
$$
then there exists a $\tau \in (-r^2,-\a r^2)$, such that
$$
\mu(\tau)\geq (\fr 12-\a)(1-\a)^{-1}
mes(S_{r})mes(K_{r}),\leqno(3.6)
$$
we have by noticing $v=0$ when $u\geq 1,$
$$
\int_{K_{r}}\int_{S_{\beta r}} v(\tau,x',\overline x)d \overline x
dx'\leq \fr 12(1-\a)^{-1}mes(S_{r})mes(K_{r})\ln(h^{-\fr 9
8}).\leqno(3.7)
$$
Now we choose $\varepsilon={\fr{1}{2\la}}\,$ and $\a$  (near zero)
and $\beta$ (near one), so that
$$
\fr{1}{4\beta^{2Q}}+\fr{1}{2\beta ^{2Q}(1-\a)}\leq \fr 4
5.\leqno(3.8)
$$
By (3.2), (3.5), (3.7) and (3.8), and note the last term in (3.5)
can be controlled by $C(B,\la,N)(1-\beta)^{-2}\beta^{-Q}|K_{\beta
r}| | S_{\beta r}|,$ we deduce
$$\ba{lll}
\int_{K_{\beta r}}\int_{S_{\beta r}} v(t,x',\overline x)d \overline
x dx'\\ \\ \leq [2C(1-\beta)^{-2}\beta^{-Q} +\fr 45\ln(h^{-\fr 9
8})]mes(K_{\beta r}\times S_{\beta r}).\ea\leqno(3.9)
$$
When $(x', \bar{x})\notin {\cal N}_{t,h},$ $u\geq h$, we have
$$\ln(\fr 1 {2h})\leq \ln^+(\fr{1}{h+h^{\fr 9 8}})\leq v,$$
then
$$\ln(\fr 1
{2h})mes(K_{\beta r}\times S_{\beta r}\setminus {\cal N}_{t,h})\leq
\int_{K_{\beta r}}\int_{S_{\beta r}} v(t,x',\overline x)d \overline
x dx'.$$
Since
$$
\fr{C+{\fr 45}\ln(h^{-\fr 98})}{\ln(h^{-1})}\lra \fr
9{10},\qquad\hbox{as} \quad h\ra 0,
$$
then there exists constant $h_1$ such that for $0<h<h_1$ and $t
\in(-\a r^2,0)$
$$
mes(K_{\beta r}\times S_{\beta r}\setminus {\cal N}_{t,h})\leq \fr
{10}{11}mes(K_{\beta r}\times S_{\beta r}).
$$
Then we proved our lemma.
\begin{cor}
Under the assumptions of Lemma 3.1, we can choose $\theta$,
$0<\theta< \a$ and $\theta< \beta$ small enough so that
$$
mes\{{\cal B}^-_{\beta r} \setminus {\cal B}^-_{\theta r} \cap
\{(t,x)|\quad u\geq h\}\}\geq C_0(\a,\beta,\Lambda) mes \{{\cal
B}^-_{\beta r}\},
$$
where $0<C_0(\a,\beta,\Lambda)<1$.
\end{cor}

Let $\chi(s)$ be a smooth function given by
$$\ba{ll}
\chi(s)=1 \qquad if \quad s\leq {\theta^{\fr 1 {2Q}}} r,\\
\chi(s)=0 \qquad if \quad s> r, \ea
$$
where ${\theta^{\fr 1 {2Q}}}<\fr {1}{2}$ is a constant. Moreover, we
assume that
$$
0\leq -\chi'(s) \leq \fr{2}{(1 -{\theta^{\fr 1 {2Q}}})r},
$$
and $\chi'(s)<0$, if ${\theta^{\fr 1 {2Q}}} r<s< r$. Also for any
$\beta_1, \beta_2,$ with $\theta^{\fr 1 {2Q}}<\beta_1<\beta_2<1,$ we
have $$|\chi'(s)|\geq C(\beta_1,\beta_2)>0,$$ if $\beta_1r\leq s\leq
\beta_2r.$

For $x\in R^N,$ $t<0$, we set $${\mathcal Q} =\{(x',\bar{x},t)|
-r^2\leq t< 0,\, x'\in K_{\fr r \theta},
\,|x_j|\leq\fr{r^{\a_j}}{\theta}, j=m_0+1,\cdots, N\},$$
$$
\phi_0(x,t)=\chi([\theta^2|t|^Q\langle {\mathcal
C}^{-1}(|t|)e^{tB^T}x,e^{tB^T}x\rangle+\theta^2\sum_{i=m_0+1}^{N}
\fr {x_i^2}{r^{2\a_i-2Q}}-C_1 t r^{2Q-2} ]^{\fr {1} {2Q}}),
$$
$$
\phi_1(x)=\chi(\theta |x'|),
$$
$$
\phi(t,x)=\phi_0(t,x)\phi_1(x),\leqno(3.10)
$$
where $C_1>1$ is chosen so that
$$
\ba{lllll}  C_1 r^{2Q-2}&\geq 2\theta^2|t|^Q|\langle x,B e^{t
B}{\mathcal C}^{-1}(|t|)e^{t B^T}x\rangle |\\
\\&+\theta^2|t|^{Q}\langle {\mathcal C}^{-1}(|t|)e^{t B^T}x,A_0
{\mathcal C}^{-1}(|t|)e^{t B^T}x\rangle\\
\\ &+\theta^2|\sum_{i=1}^{N}\sum_{j>m_0}2x_i b_{ij}x_j
r^{2Q-2\alpha_j}|, \ea
$$
for all $z\in {\mathcal Q}$.

In the following discussion, $a\approx b$ means
$$C(B,\la,N)^{-1}a\leq b\leq C(B,\la,N)a.$$

With the notations given in section 2, for $ s>0$, or $t<0$, we have
following properties:
$${\mathcal C}'(s)=A_0-B^T {\mathcal C}(s)-{\mathcal C}(s)B, \leqno(a)$$
$$Y\langle {\mathcal C}^{-1}(|t|)x,x\rangle=4\langle
x,B{\mathcal C}^{-1}(|t|)x\rangle- \langle {\mathcal
C}^{-1}(|t|)x,A_0 {\mathcal C}^{-1}(|t|)x\rangle,\leqno(b) $$
 $$\ba{lll}Y\langle {\mathcal C}^{-1}(|t|)e^{t B^T}x,e^{t
B^T}x\rangle=&2\langle x,B e^{t B}{\mathcal C}^{-1}(|t|)e^{t
B^T}x\rangle \\\\ {}&- \langle {\mathcal C}^{-1}(|t|)e^{t B^T}x,A_0
{\mathcal C}^{-1}(|t|)e^{t B^T}x\rangle;\ea\leqno(c)
$$
(d) moreover, if $|t|$ is small enough, then
$$
\hspace*{12pt}\langle {\mathcal C}^{-1}(|t|)e^{t B^T}x,e^{t
B^T}x\rangle\approx |D_{|t|^{-\fr 1 2}}x|^2,\leqno(d.1)$$
$$
\langle {\mathcal C}^{-1}(|t|)e^{t B^T}B^Tx,e^{t B^T}x\rangle\leq C
|t|^{-1} |D_{|t|^{-\fr 1 2}}x|^2,\leqno(d.2)
$$
$$
\langle A_0 {\mathcal C}^{-1}(|t|)e^{t B^T}x,{\mathcal
C}^{-1}(|t|)e^{t B^T}x\rangle\leq C |t|^{-1} |D_{|t|^{-\fr 1
2}}x|^2, \leqno(d.3)
$$
where C depends on B, $\la$, and N.

The property (a) can be checked by the definition,  in fact,
$${\mathcal
C}(s)=\int_0^s E(t)A_0 E^T(t)\rm{d}t,$$
then
$$
{\mathcal C}'(s)=E(s)A_0 E^T(s),
$$
$$
{\mathcal C}''(s)=E(s)(-B^T)A_0 E^T(s)+E(s)A_0 E^T(s)(-B)=-B^T
{\mathcal C}'(s)-{\mathcal C}'(s)B,
$$
integrating from 0 to s,\, we have
$$
{\mathcal C}'(s)=A_0-B^T {\mathcal C}(s)-{\mathcal C}(s)B.
$$
To prove (b), we calculate
\begin{eqnarray*}
\lefteqn{Y\langle {\mathcal C}^{-1}(|t|)x,x\rangle}\,\\
&&=[\langle x,B D\rangle-\partial_t]\langle {\mathcal C}^{-1}(|t|)x,x\rangle\,\\
&&=2\langle x,B{\mathcal C}^{-1}(|t|)x\rangle+\langle \partial_{|t|}{\mathcal C}^{-1}(|t|)x,x\rangle\,\\
&&=2\langle x,B{\mathcal C}^{-1}(|t|)x\rangle-\langle {\mathcal
C}^{-1}(|t|)
\partial_{|t|}{\mathcal C}(|t|){\mathcal C}^{-1}(|t|)x,x\rangle\,\\
&&=4\langle x,B{\mathcal C}^{-1}(|t|)x\rangle- \langle {\mathcal
C}^{-1}(|t|)x,A_0 {\mathcal C}^{-1}(|t|)x\rangle.
\end{eqnarray*}
The proof of (c) is the same as (b). \\
Applying (2.7) and (2.5),
\begin{eqnarray*}
\langle {\mathcal C}^{-1}(|t|)e^{t B^T}x,e^{t B^T}x\rangle &\approx
&\langle {\mathcal C}_0^{-1}(|t|)e^{t B^T}x,e^{t B^T}x\rangle\\ &=&
\langle {\mathcal C}_0^{-1}(1)D_{|t|^{-\fr 1
2}}e^{t\,B^T}x,D_{|t|^{-\fr 1 2}}e^{t\,B^T}x\rangle\\ &\approx
&\|e^{\tilde{B}}D_{|t|^{-\fr 1 2}}x\|\approx |D_{|t|^{-\fr 1 2}}x|^2
\end{eqnarray*}
where $D_{|t|^{-\fr 1 2}}B^T=|t|^{-1}\tilde{B}D_{|t|^{-\fr 1 2}},\,$
$D_{|t|^{-\fr 1 2}}e^{tB^T}=e^{\tilde{B}} \, D_{|t|^{-\fr 1 2}}\,$
and $\tilde{B}$\, has the form
$$\left(
\begin{array}{ccccc}
|t|B_{0,0}^T  &|t|^2 B_{1,0}^T & \cdots & \cdots & |t|^{d+1}B_{d,0}^T  \\
B_1^T &  |t|B_{1,1}^T   &\cdots &\cdots & |t|^{d}B_{d,1}^T \\
0 &B_2^T& \ddots &\cdots &\vdots \\
\vdots & \ddots &\ddots &\ddots & \vdots \\
0 &  \cdots  & 0 &  B_d^T &|t|B_{d,d}^T \\
\end{array}
\right)
$$
B is given by
$$\left(
\begin{array}{ccccc}
B_{0,0}&B_1&0&\cdots&0\\
B_{1,0}&B_{1,1}&B_2&\cdots&0\\
\vdots &\vdots &\vdots &\ddots &\vdots \\
B_{d-1,0}&B_{d-1,1}&B_{d-1,2}&\cdots&B_d\\
B_{d,0}&B_{d,1}&B_{d,2}&\cdots&B_{d,d}\\
\end{array}
\right)
$$
then we obtain (d.1).\\
For any  $x\in R^N,$ by the Young inequality, (2.5) and (2.7), we
have
$$
\ba{llllllll}
&&\langle {\mathcal C}^{-1}(|t|)B^Te^{t B^T}x,e^{t B^T}x\rangle\\
\\ &\leq & \varepsilon\langle {\mathcal C}^{-1}(|t|)B^Te^{t B^T}x,B^Te^{t
B^T}x\rangle+\fr{1}{4\varepsilon}\langle {\mathcal C}^{-1}(|t|)e^{t
B^T}x,e^{t
B^T}x\rangle \\ \\
&\leq&2(\varepsilon\langle {\mathcal C}_0^{-1}(|t|)B^Te^{t
B^T}x,B^Te^{t B^T}x\rangle+\fr{1}{4\varepsilon}\langle {\mathcal
C}_0^{-1}(|t|)e^{t B^T}x,e^{t
B^T}x\rangle) \\ \\
&\leq& C(B,\la,N)(\varepsilon |t|^{-2}|D_{|t|^{-\fr 1
2}}x|^2+\fr{1}{4\varepsilon}|D_{|t|^{-\fr 1 2}}x|^2)\\ \\
&=&C(B,\la,N) |t|^{-1}|D_{|t|^{-\fr 1 2}}x|^2\quad
(\varepsilon=|t|), \ea
$$
which is (d.2).

Let $y=e^{t B^T}x$, by (2.5), (2.6) and (2.7), we have
$$
\ba{llllllll} & &\langle y, {\mathcal C}^{-1}(|t|)A_0
{\mathcal C}^{-1}(|t|)y\rangle\\ \\
&=&\langle y, {\mathcal C}^{-1}(|t|)A_0 A_0
{\mathcal C}^{-1}(|t|)y\rangle\\ \\
&=&|t|^{-1}\langle y, {\mathcal C}^{-1}(|t|)A_0D_{|t|^{\fr 1
2}}D_{|t|^{\fr 1 2}} A_0
{\mathcal C}^{-1}(|t|)y\rangle\\ \\
&\leq&|t|^{-1}\langle y, {\mathcal C}^{-1}(|t|)D_{|t|^{\fr 1
2}}D_{|t|^{\fr 1 2}} {\mathcal C}^{-1}(|t|)y\rangle \ea
$$
consequently,
$$ \ba{llllllll}& &|t|^{-1}\langle y, {\mathcal C}^{-1}(|t|)D_{|t|^{\fr 1
2}}D_{|t|^{\fr 1 2}} {\mathcal C}^{-1}(|t|)y\rangle\\ \\
&\approx&|t|^{-1}\langle y, {\mathcal C}^{-1}(|t|)D_{|t|^{\fr 1
2}}{\mathcal C}_0(1)D_{|t|^{\fr 1 2}}
{\mathcal C}^{-1}(|t|)y\rangle\\ \\
&=&|t|^{-1}\langle y, {\mathcal C}^{-1}(|t|){\mathcal C}_0(|t|)
{\mathcal C}^{-1}(|t|)y\rangle\\ \\
&\approx&|t|^{-1}\langle y, {\mathcal C}^{-1}(|t|){\mathcal C}(|t|)
{\mathcal C}^{-1}(|t|)y\rangle\\ \\
&=&|t|^{-1}\langle y,
{\mathcal C}^{-1}(|t|)y\rangle\\ \\
&\approx&|t|^{-1}\langle y,
{\mathcal C}_0^{-1}(|t|)y\rangle\\ \\
&\approx&|t|^{-1}|D_{|t|^{-\fr 1 2}}y|^2\\ \\
&\approx&|t|^{-1}|D_{|t|^{-\fr 1 2}}x|^2 \ea
$$
and we obtain the proof of (d.3).
\begin{rmk} By the definition of $\phi$ and the above arguments,
it is easy to check that, for $\theta$,  $r$ small  and $t\leq 0$\\
(1) $\phi(z)\equiv 1,$ in ${\cal B}^-_{\theta r}$,\\
(2) $\rm{supp}\phi\subset {\mathcal Q}$,\\
(3) there exists $\a_1>0,$ which depends on $C_1,$ such that
$$\{(x,t)| -\a_1r^2\leq t < 0, x'\in K_r, \bar{x}\in S_{\beta
r}\}\subseteq \rm{supp}\phi, $$ (4) $0<\phi_0(z)<1,$ for $z\in
\{(x,t)| -\a_1r^2\leq t \leq -\theta r^2, x'\in K_r, \bar{x}\in
S_{\beta r}\}$.
\end{rmk}
 \begin{lemma} Under the above notations, we have \\
(e) $Y \phi_0(z)\leq 0, \quad \rm{for}\quad z\in {\mathcal Q}$;\\
(f) $$|\int_{\mathcal {Q}}\phi_1|\langle A_0 D\phi_0,D
\Gamma_1(z,\cdot)\rangle|-\int_{\mathcal {Q}}\phi_1\langle A_0
D\phi_0,D \Gamma_1(0,\cdot)\rangle|\leq C_6 \theta^2,
$$\, \rm{for} $z\in{\cal B}^-_{\tilde{\theta} r}$, where $C_6$ is
dependant on $B$, $\la$, $N$ and $\tilde{\theta}$ depends on
$\theta$.
\end{lemma}
{\it Proof:} \\Let
$$
[\theta^2|t|^Q\langle {\mathcal
C}^{-1}(|t|)e^{tB^T}x,e^{tB^T}x\rangle+\theta^2\sum_{i=m_0+1}^{N}
\fr {x_i^2}{r^{2\a_i-2Q}}-C_1 t r^{2Q-2}]
$$
be denoted by $[\cdots]$. Then
$$\ba{llllllllll} Y \phi_0 &=
\chi'([\cdots]^{\fr{1}{2Q}})\fr{1}{2Q}[\cdots]^{\fr{1}{2Q}-1}
[\theta^2|t|^QY\langle {\mathcal C}^{-1}(|t|)e^{t B^T}x,e^{t
B^T}x\rangle\,\\ \\
 & + Q\theta^2|t|^{Q-1}\langle {\mathcal C}^{-1}(|t|)e^{t B^T}x,e^{t B^T}x\rangle+C_1
r^{2Q-2}\\ \\
&+\theta^2\sum_{i=1}^{N}\sum_{j>m_0}(2x_i b_{ij}x_j
r^{2Q-2\alpha_j})]
\,\\ \\
&=\chi'([\cdots]^{\fr{1}{2Q}})\fr{1}{2Q}[\cdots]^{\fr{1}{2Q}-1}
[\theta^2|t|^Q(2\langle x,B e^{t B}{\mathcal C}^{-1}(|t|)e^{t
B^T}x\rangle \,\\ \\
& -\langle {\mathcal C}^{-1}(|t|)e^{t B^T}x,A_0 {\mathcal
C}^{-1}(|t|)e^{t B^T}x\rangle)\, +Q\theta^2|t|^{Q-1}\langle
{\mathcal C}^{-1}(|t|)e^{t B^T}x,e^{t
B^T}x\rangle \\ \\
&  +C_1 r^{2Q-2}+\theta^2\sum_{i=1}^{N}\sum_{j>m_0}(2x_i b_{ij}x_j
r^{2Q-2\alpha_j})]. \, \ea
$$
We choose $C_1>1,$ such that
$$
\ba{lll}  C_1 r^{2Q-2}\geq &\theta^2|t|^Q(2|\langle x,B e^{t
B}{\mathcal C}^{-1}(|t|)e^{t B^T}x\rangle |\\ \\&+\langle {\mathcal
C}^{-1}(|t|)e^{t B^T}x,A_0
{\mathcal C}^{-1}(|t|)e^{t B^T}x\rangle)\,\\
\\ &+\theta^2|\sum_{i=1}^{N}\sum_{j>m_0}2x_i b_{ij}x_j
r^{2Q-2\alpha_j}|, \ea
$$
by the above (d),
$$
\theta^2|t|^Q|\langle x,B e^{t B}{\mathcal C}^{-1}(|t|)e^{t
B^T}x\rangle |\leq C\theta^2|t|^{Q-1}|D_{|t|^{-\fr 1 2}}x|^2\leq C
r^{2Q-2},
$$
for all $z\in {\mathcal Q}.$

Similar results holds for $\theta^2|t|^Q \langle {\mathcal
C}^{-1}(|t|)e^{t B^T}x,A_0 {\mathcal C}^{-1}(|t|)e^{t B^T}x\rangle.$
For the term $x_i b_{ij}x_j r^{2Q-2\alpha_j}$, then either $\a_i\geq
\a_j$ or $\a_j=\a_i+2$, we also obtain
$$ \theta^2|\sum_{i=1}^{N}\sum_{j>m_0}2x_i b_{ij}x_j
r^{2Q-2\alpha_j}| \leq C(B,\la,N) r^{2Q-2}.
$$
Thus $C_1(B,\la,N)$ is well defined,
then $Y \phi_0(z)\leq 0$ ($z\in {\mathcal Q}$) holds. \\

For the proof of (f), let $g(z)=\int_{\mathcal {Q}}\phi_1|\langle
A_0 D\phi_0,D \Gamma_1(z,\cdot)\rangle|(\zeta)$, then $g(z)$ is
smooth and $g(z)\leq g(0)+C(B,\la,N)|z|$. In fact,
$$
\ba{llllllll}
&&g(0)\\\\&=&\int_{\mathcal {Q}}\phi_1\langle A_0 D\phi_0,D \Gamma_1(0,\cdot)\rangle\\
\\&=&\int_{\mathcal {Q}}\phi_1\chi'([\cdots]^{\fr{1}{2Q}})\fr{[\cdots]^{\fr{1}{2Q}-1}}{2Q}
\theta^2|\tau|^Q\langle D_{m_0} \langle {\mathcal
C}^{-1}(|\tau|)e^{\tau B^T}\xi,e^{t
B^T}\xi\rangle,\,D_{m_0}^{(\zeta)}\Gamma_1(0,\cdot)\rangle \\ \\
&=&\int_{\mathcal
{Q}}\phi_1|\chi'([\cdots]^{\fr{1}{2Q}})|\fr{\theta^2|\tau|^Q\Gamma_1(0,\zeta)}{2Q[\cdots]^{1-\fr{1}{2Q}}}
\langle e^{\tau B}{\mathcal C}^{-1}(|\tau|)e^{\tau B^T}\xi,A_0
e^{\tau B}{\mathcal C}^{-1}(|\tau|)e^{\tau B^T}\xi\rangle. \ea
$$
We choose a domain $\mathcal {D}$ as in Remark 3.1, and ${\mathcal
D}=\{(x,t)| -\a_1r^2\leq t \leq -{\fr{\a_1}2} r^2, x'\in K_r,
\bar{x}\in S_{\beta r}\},$ then by choosing small $\theta$ we get,
$0<\phi_0<1,$ $\chi'([\cdots]^{\fr{1}{2Q}})\approx r^{-1},$
$\phi_1\equiv1,$ $[\cdots]\approx r^{2Q},$ and
$\Gamma_1(0,\zeta)\approx |\tau|^{-\fr Q 2}$ when $\zeta\in \mathcal
{D}.$ Hence $$ g(0)\geq C(B,\la,N)\theta^2r^{-Q}\int_{\mathcal
{D}}\langle e^{\tau B}{\mathcal C}^{-1}(|\tau|)e^{\tau B^T}\xi,A_0
e^{\tau B}{\mathcal C}^{-1}(|\tau|)e^{\tau B^T}\xi\rangle.$$ By
$D_{|\tau|^{-\fr 1 2}}e^{\tau B^T}=e^{\tilde{B}} \, D_{|\tau|^{-\fr
1 2}}$ in $(d.1)$, and $D_{|\tau|^{\fr 1 2}}{\mathcal
C}^{-1}(|\tau|)D_{|\tau|^{\fr 1 2}}$ which is positive and whose
eigenvalues can be controlled by constants from (2.5) and (2.7),
then
$$\ba{llllll}
&&r^2\langle e^{\tau B}{\mathcal C}^{-1}(|\tau|)e^{\tau B^T}\xi,A_0
e^{\tau B}{\mathcal C}^{-1}(|\tau|)e^{\tau B^T}\xi\rangle
\\\\&=&r^2|\tau|^{-1}\langle e^{\tilde{B}^T}D_{|\tau|^{\fr 1 2}}{\mathcal
C}^{-1}(|\tau|)D_{|\tau|^{\fr 1 2}}e^{\tilde{B}}D_{|\tau|^{-\fr 1
2}}\xi, A_0e^{\tilde{B}^T}D_{|\tau|^{\fr 1 2}}{\mathcal
C}^{-1}(|\tau|)D_{|\tau|^{\fr 1 2}}e^{\tilde{B}}D_{|\tau|^{-\fr 1
2}}\xi\rangle, \ea
$$
which is positive and not dependent on $r$ except zero measurable
set, hence we get $g(0)\geq C_6\theta^2>0$ with $C_6$ as a constant
dependant on B, $\la$, N. We can choose $\tilde{\theta}$ small,
$0<\tilde{\theta}<\theta,$ such that $g(z)\leq g(0)+\fr12
C_6\theta^2$ for $z\in {\cal B}^-_{\tilde{\theta}
r}.$\\

We now have the following ${\rm Poincar\acute{e}}$'s type
inequality.
\begin{lemma}
Let $w$ be a non-negative weak sub-solution of (1.2) in ${\cal
B}_1^-$. Then there exists a constant $C$, only depends on $B,$
$\la$ and $N$, such that for $r<\theta<1$
$$
\int_{{\cal B}^-_{\theta r}}(w(z)-I_0)_+^2\leq C\theta^2
r^2\int_{{\cal B}^-_{\fr r {\theta}}}|D_{m_0}w|^2, \leqno(3.11)
$$
where $I_0$ is given by
$$
I_0=max_{{\cal B}^-_{\tilde{\theta} r}}[I_1(z)+C_2(z)],\leqno(3.12)
$$
and
$$
I_1(z)=\int_{{\cal B}^-_{\fr r {\theta}}} [\langle
{\phi}_1A_0D{\phi}_0,D\Gamma_1(z,\cdot)\rangle
w-\Gamma_1(z,\cdot)wY\phi](\zeta)d\zeta,\leqno(3.13)
$$
$$
C_2(z)=\int_{{\cal B}^-_{\fr r {\theta}}} [\langle
{\phi}_0A_0D{\phi}_1,D\Gamma_1(z,\cdot)\rangle w](\zeta)d\zeta,
$$
where $\Gamma_1$ is the fundamental solution,  and $\phi$ is given
by (3.10).
\end{lemma}
 {\it Proof:} We represent $w$ in terms of the fundamental
solution of $\Gamma_1$. For $z \in {\cal B}^-_{\theta r}$, we have
$$\ba{llll}
w(z)&=\int_{{\cal B}^-_{\fr r {\theta}}}  [\langle
A_0D(w\phi),D\Gamma_1(z,\cdot)\rangle
-\Gamma_1(z,\cdot)Y(w\phi)](\zeta)d\zeta \\ \\&=
I_1(z)+I_2(z)+I_3(z)+C_2(z),\ea\leqno(3.14)
$$
where $I_1(z)$ and $C_2(z)$ are given by (3.13) and
$$
I_2(z)=\int_{{\cal B}^-_{\fr r {\theta}}} [\langle
(A_0-A)Dw,D\Gamma_1(z,\cdot)\rangle\phi-\Gamma_1(z,\cdot)\langle
ADw,D\phi\rangle](\zeta)d\zeta,
$$
$$
I_3(z)=\int_{{\cal B}^-_{\fr r {\theta}}} [\langle
ADw,D(\Gamma_1(z,\cdot)\phi)\rangle-\Gamma_1(z,\cdot)\phi
Yw](\zeta)d\zeta.
$$
From our assumption, $w$ is a weak sub-solution of (1.2), and $\phi$
is a test function of this semi-cylinder. In fact, we let
$$
\tilde{\chi}(\tau)=\left\{
\begin{array}{lll} 1\quad &\tau\leq
0,\\1-n\tau\quad &0\leq \tau \leq 1/n,\\ 0\quad &\tau\geq
1/n.\end{array}\right.$$ Then
$\tilde{\chi}(\tau)\phi\Gamma_1(z,\cdot)$ can be a test function
(see [2]). Let $n\rightarrow \infty$, we obtain
$\phi\Gamma_1(z,\cdot)$ as a legitimate test function, and
$I_3(z)\leq 0$. Then in ${\cal B}^-_{\theta r}$,
$$
0\leq (w(z)-I_0)_+\leq I_2(z)=I_{21}+I_{22}.
$$
By Corollary 2.1 we have
$$
||I_{21}||_{L^2({\cal B}^-_{\theta r})}\leq C(\la,N)\theta
r||I_{21}||_{L^{2+\fr 4 Q}({\cal B}^-_{\theta r})}\leq C(B, \la, N)
\theta r||D_{m_0}w||_{L^2({\cal B}^-_{{\fr r
{\theta}}})}.\leqno(3.15)
$$
Similarly for $I_{22},$
$$
||I_{22}||_{L^2({\cal B}^-_{\theta r})}\leq |{\cal B}^-_{\theta
r}|^{\fr 12-\fr{Q-2}{2Q+4}} ||I_{22}||_{L^{2\tilde{k}}({\cal
B}^-_{\theta r})}\leq C(B,\la,N) \theta^2 r^2||D_{m_0}w
D_{m_0}\phi||_{L^2({\cal B}^-_{{\fr r {\theta}}})},
$$
where $D_{m_0}\phi=\phi_0 D_{m_0}\phi_1+\phi_1D_{m_0}\phi_0.$
$$
|\phi_0 D_{m_0}\phi_1|=|\phi_0\chi'(\theta|\xi'|)\theta
D_{m_0}(|\xi'|)|\leq \fr{\theta}{r},
$$
and
$$
\ba{llllll} |\phi_1D_{m_0}\phi_0|&\leq& 2
\phi_1|\chi'([\cdots]^{\fr{1}{2Q}})|\fr{1}{2Q}[\cdots]^{\fr{1}{2Q}-1}\theta^2
|\tau|^Q|A_0e^{\tau B}{\mathcal C}^{-1}(|t|)e^{\tau B^T}\xi|\\\\
&\leq& C(B,\la,N)r^{-1}(\theta
r^{2Q})^{\fr{1}{2Q}-1}\theta^2|\tau|^{Q-\fr 12}|D_{|\tau|^{-\fr 1
2}}\xi|\\ \\
&\leq& C(B,\la,N)r^{-1}(\theta
r^{2Q})^{\fr{1}{2Q}-1}\theta^2 r^{2Q-1}\theta^{-1}\\ \\
&\leq& C(B,\la,N)\theta^\fr{1}{2Q} r^{-1}, \ea
$$
thus $$||I_{22}||_{L^2({\cal B}^-_{\theta r})}\leq C(B,\la,N)
\theta^2 r||D_{m_0}w||_{L^2({\cal B}^-_{{\fr r {\theta}}})}.$$\\
Then we proved our lemma.

Now we apply Lemma 3.3 to the function
$$
w= \ln^+\fr{h}{u+h^{\fr 98}}.
$$
If $u$ is a weak solution of (1.2), obviously $w$ is a weak
sub-solution. We estimate the value of $I_0$ given by (3.12) and
(3.13) in Lemma 3.3.
\begin{lemma}
Under the assumptions of Lemma 3.3, there exist constants
$\lambda_0$, $r_0$ and $r_0<\theta$. $\la_0$ only depends on
constants $\a$,
 $\beta$, $\lambda$, $B$, $N$, and $\phi$, $0<\lambda_0<1$, such
 that for $r<r_0$
$$
|I_0|\leq \lambda_0 \ln(h^{-\fr 1 8}).\leqno(3.16)
$$
\end{lemma}
{\it Proof:} We first come to estimate $C_2(z)$ and often denote
$x=(x',\bar{x},t)$, and $\zeta=(\xi',\bar{\xi},\tau)$. Note
supp$\phi\in {\mathcal Q}$, and $z \in {\cal B}^-_{\theta r}$, then
$$
\ba{llllllllll} & & |C_2(z)|\\\\&=& |\int_{{\cal B}^-_{\fr r
{\theta}}} [\langle {\phi}_0A_0D{\phi}_1,D\Gamma_1(z,\cdot)\rangle
w](\zeta)d\zeta|
\\
&\leq & C(B,\la,N) \ln (h^{-\fr 1 8})\fr{2m_0\theta}{(1-\theta^{\fr
1 {2Q}})r} \sup_{\theta|\xi'|\geq\theta^{\fr 1 {2Q}}r
}||\zeta^{-1}\circ z||^{-Q-1}\cdot \theta^{-N}|r|^{Q+2}
\\
&\leq & C(B,\la,N) \ln (h^{-\fr 1 8})\fr{2m_0}{(1-\theta^{\fr 1
{2Q}})r}\theta |\theta^{{\fr 1 {2Q}}-1}r-\theta r|^{-Q-1}\cdot\,
\theta^{-N}|r|^{Q+2}\\
&\leq & C(B,\la,N) \theta^{Q+{\fr 3 2}-N-{\fr 1 {2Q}}}\ln( h^{-\fr 1 8})\\
&=& C_3 \theta^{\a_0} \ln (h^{-\fr 1 8})\, \ea \leqno(3.17)
$$
where $\a_0=Q+{\fr 3 2}-N-{\fr 1 {2Q}} >0$ and similarly
$$
\ba{lllllllll} & &|\int_{{\cal B}^-_{\fr r {\theta}}} [-\phi_0
Y\phi_1\Gamma_1(z,\cdot) w](\zeta)d\zeta| \\&\leq & |\int_{{\cal
B}^-_{\fr r {\theta}}} [-\phi_0 \chi'(\theta|\xi'|)\theta
\sum_{i=1}^{N}\sum_{j=1}^{m_0} \xi_i
b_{ij}\xi_j/|\xi'|\,\Gamma(z,\cdot) w](\zeta)d\zeta|
\\
&\leq & C(B,\la,N) |\theta^{{\fr 1 {2Q}}-1}r-\theta
r|^{-Q}\theta^{-N}|r|^{Q+2} \ln \,(h^{-\fr 1 8})\\ &\leq &C(B,\la,N)
\theta^{Q-N-{\fr 1 2}} r^2 \ln (h^{-\fr 1
8})\,\\
&=&C_4 \theta^{Q-N-{\fr 1 2}} r^2 \ln (h^{-\fr 1 8})\,
\\&\leq &C_4 \theta^{\tilde{\a_0}} \ln (h^{-\fr 1 8})\, \ea\leqno (3.18)
$$
where $\tilde{\a_0}=Q+{\fr 3 2}-N>0$, if $r<\theta$.

Now we let $w\equiv 1$ then (3.14) gives, for $z \in {\cal
B}^-_{\theta r}$,
$$
\ba{llll} 1&=&\int_{{\cal B}^-_{\fr r {\theta}}} [\langle
\phi_1A_0D\phi_0,D\Gamma_1(z,\cdot)\rangle-
\phi_1\Gamma_1(z,\cdot)Y\phi_0](\zeta)d\zeta\\
\\&&-\int_{{\cal B}^-_{\fr r {\theta}}}\phi_0\Gamma_1(z,\cdot)Y\phi_1(\zeta)d\zeta
+C_2(z)|_{w=1},\ea\ \leqno(3.19)
$$
where $\phi$ is given by (3.10). By Lemma 3.2, for $z\in
B^-_{\tilde{\theta} r}$,
$$
-\phi_1\Gamma_1(z,\cdot)Y\phi_0 \geq 0.\leqno(3.20)
$$
$$
g(z)=\int_{\mathcal Q}|\langle \phi_1
A_0D\phi_0,D\Gamma_1(z,\cdot)\rangle|\leq g(0)+\fr12
C_6\theta^2.\leqno(3.21)
$$
We only need to prove $-\phi_1\Gamma_1(z,\cdot)Y\phi_1$ has a
positive lower bound in a domain which $w$ vanishes, and this bound
independent of $\,r$ and small $\theta$. So we can find a $\la_0,$
$0<\la_0<1$, such that this lemma holds and $\la_0$ is independent
of $r$ and small $\theta.$ We observe that the support of $\chi'(s)$
is in the region ${\theta^{\fr 1 {2Q}}}r<s< r$, thus for some
${\beta}'< 1$ (we choose $\beta'$ near one), the set ${\cal
B}^-_{{\beta}' r} \setminus {\cal B}^-_{{\sqrt\theta} r}$ with
$|t|>\theta r^2/{C_1}$ is contained in the support of $\phi'$ . Then
we can prove that the integral of (3.20) on the domain ${\cal
B}^-_{{\beta}' r} \setminus {\cal B}^-_{{\sqrt\theta} r}$ with
$|t|>\theta r^2/{C_1}$ is lower bounded by a positive constant.

For $z\in B^-_{\theta r}$, $0<\a_1\leq \a$ and set
$$\zeta\in Z=\{(\xi,\tau)| -\a_1 r^2\leq \tau\leq -\fr {\a_1}
{2}r^2,\,x'\in K_r, \,\bar{x}\in S_{\beta r},\,w(\xi,\tau)=0\},$$
then $|Z|=C(\a_1, \la,N)r^{Q+2}$ by Lemma 3.1 and Corollary 3.1. We
note that $w(\zeta)=0,$ $\phi_1(\zeta)=1$,
$|\chi'([\cdots]^{\fr{1}{2Q}})|\geq C(\a_1,B,\la,N)>0$ when
$\zeta\in Z $ and $\theta$ is small, then
$$
\ba{llllllllllll} \int_Z [-\phi_1\Gamma_1(z,\cdot)Y\phi_0](\zeta)\,d \zeta \\
\\= -\int_Z \phi_1\Gamma_1(z,\cdot)
\chi'([\cdots]^{\fr{1}{2Q}})\fr{1}{2Q}[\cdots]^{\fr{1}{2Q}-1}
[\theta^2|\tau|^Q(2\langle \xi,B e^{\tau B}{\mathcal C}^{-1}(|\tau|)e^{\tau B^T}\xi\rangle \\ \\
- \langle {\mathcal C}^{-1}(|\tau|)e^{\tau B^T}\xi,A_0 {\mathcal
C}^{-1}(|\tau|)e^{\tau B^T}\xi\rangle)
 + Q\theta^2|\tau|^{Q-1}\langle {\mathcal C}^{-1}(|\tau|)e^{\tau B^T}\xi,e^{\tau B^T}\xi\rangle \,\\ \\
 +C_1
r^{2Q-2}
+\theta^2\sum_{i=1}^{N}\sum_{j>m_0}(2\xi_i b_{ij}\xi_j r^{2Q-\alpha_j})]d\zeta\\ \\
\geq C(B,\la,\a,N)\int_Z r^{2Q-2}[r^{2Q}]^{\fr{1}{2Q}-1}
r^{-1}\Gamma_1(\zeta^{-1}\circ z;0)d\zeta \\
\\ \geq C(B,\la,\a,N)   \int_Z
r^{-2} (t-\tau)^{-{\fr Q 2}}\exp({-C|D_{{|t-\tau|}^{-\fr 12}}(x-E(t-\tau)\xi)|^2})\\ \\
\geq C(B,\la,\a,N)   \int_Z
r^{-2} (t-\tau)^{-{\fr Q 2}}\exp({-C|D_{{|\tau|}^{-\fr 12}}\xi|^2})\quad (\rm{the\,\,same\,\,as} (d.1))\\ \\
\geq C(B,\la,\a,N)\int_Z r^{-Q-2}\\\\=C(B,\la,\a,\beta,N)=C_5>0. \ea
$$
Similarly we can choose a small data, still denote $\tilde{\theta}$,
such that $$\int_{\mathcal Q}
[-\phi_1\Gamma_1(z,\cdot)Y\phi_0](\zeta)\,d \zeta\leq \int_{\mathcal
Q} [-\phi_1\Gamma_1(0,\cdot)Y\phi_0](\zeta)\,d
\zeta+\fr12C_6\theta^2,$$ for $z\in B^-_{\tilde{\theta} r}$. We can
choose a small $\theta$ ,then choose $\tilde{\theta}$ small, and
fixed them from now on, $r_0<\theta$, such that
$$\ba{llllll}
&&|I_0|\\\\&\leq& (\int_{\mathcal Q}|\langle \phi_1
A_0D\phi_0,D\Gamma_1(z,\cdot)\rangle|+
[-\phi_1\Gamma_1(z,\cdot)Y\phi_0](\zeta)\,d \zeta -C_5)\ln(h^{-{\fr
1 8
}})\\\\&&+(C_3\theta^{\a_0}+C_4\theta^{\tilde{\a_0}})\ln(h^{-{\fr 1
8 }})\\\\&\leq&
(1-C_5+C_6\theta^2+C_3\theta^{\a_0}+C_4\theta^{\tilde{\a_0}})\ln(h^{-{\fr
1 8 }})+(C_3\theta^{\a_0}+C_4\theta^{\tilde{\a_0}})\ln(h^{-{\fr 1 8
}})\\\\&\leq& \la_0 \ln(h^{-{\fr 1 8 }}).\ea
$$
Where $0<r<r_0$, $0<\la_0<1$, depends on $\a$, $\beta$, $B$, $\la$,
$N$, and $\phi$.
\begin{lemma}
Suppose that $u(x,t)\geq 0$ be a solution of equation (1.2) in
${\cal B}^-_r$ centered at $(0,0)$ and
$$
mes\{(x,t)\in {\cal B}^-_r, \quad u \geq 1\} \geq \fr 1 2 mes ({\cal
B}^-_r).
$$
Then there exist constant $\theta$ and $h_0$, $0<\theta, h_0<1$
which only depend on $B$, $\la$, $\la_0$ and $N$ such that
$$
u(x,t) \geq h_0\quad \hbox{in}\quad {\cal B}^-_{\theta r}.
$$
\end{lemma}
{\it Proof:} We consider $$w=\ln^+(\fr{h}{u+h^{\fr98}}),$$ for
$0<h<1$, to be decided. By applying Lemma 3.3 to $w$, we have $$
-\!\!\!\!\!\!\int_{{\cal B}^-_{\theta r}}( w-I_0)_+^2 \leq
C(B,\la,N)\fr{\theta r^2}{|{\cal B}^-_{\theta r}|} \int_{{\cal
B}^-_r }|D_{m_0}w|^2.
$$
Let $\tilde{u}={\fr u h}$, then $\tilde{u}$ satisfies the conditions
of Lemma 3.1. We can get similar estimates as (3.2), (3.5), (3.7)
and (3.8), hence we have
$$
\ba{lllllll} && C(B,\la,N)\fr{\theta r^2}{|{\cal B}^-_{\theta r}|}
\int_{{\cal B}^-_{ r }}|D_{m_0}w|^2\\ \\
& &\leq C(B,\la,N)\fr{\theta r^2}{|{\cal B}^-_{\theta
r}|}[C(B,\la,N)(1-\beta)^{-2}\beta^{-Q} +\fr 45\ln(h^{-\fr 18})]
mes(K_{\beta r}\times S_{\beta r})\\ \\
& &\leq C(\theta,B,N,\la) \ln(h^{-\fr 18}),
 \ea\leqno(3.22)
$$
where $\theta$ has been chosen. By Lemma 2.2, there exists a
constant, still denoted by $\theta$, such that for $z \in {\cal
B}^-_{\theta r}$,
$$
w-I_0\leq C(B,\la,N) (\ln(h^{-\fr 18}))^{\fr 12} .\leqno(3.23)
$$
Therefore we may choose $h_0$ small enough, so that
$$
C (\ln (\fr {1}{h_0^{\fr 18}}))^{\fr 12}\leq  \ln (\fr
{1}{2h_0^{\fr 18}})-\lambda_0\ln (\fr {1}{h_0^{\fr 18}}).
$$
Then (3.16) and (3.23) implies
$$
\max_{{\cal B}^-_{\theta r}}\fr{h_0}{u+h_0^{\fr 98}}\leq \fr
{1}{2h_0^{\fr 18}},
$$
which implies $\min_{{\cal B}^-_{\theta r}}u\geq h_0^{\fr 98}$, then
we finished the proof of this Lemma.\\

{\bf Proof of Theorem 1.1.} We may assume that $M=\max_{{\cal
B}^-_{r}}(+u)=\max_{{\cal B}^-_{r}}(-u)$, otherwise we replace $u$
by $u-c$, since $u$ is bounded locally. Then either $1+\fr u M$ or
$1-\fr u M$ satisfies the assumption of Lemma 3.5, and we suppose
$1+\fr u M$ does it, thus Lemma 3.5 implies existing $h_0>0$ such
than $\inf_{{\cal B}^-_{\theta r}}(1+\fr u M)\geq h_0,\, $ i.e.
$u\geq M(h_0-1)$,\,then
$$
Osc_{{\cal B}^-_{\theta r}}u\leq M-M(h_0-1)\leq
(1-\fr{h_0}{2})Osc_{{\cal B}^-_{r}}u,
$$
which implies the $C^{\a}$ regularity of $u$ near point $(0,0)$ by
the standard iteration arguments. By the left invariant
translation group action, we know that $u$ is $C^{\a}$ in the
interior.

{\bf Acknowledgments.} We would like to thank the referees for
valuable comments and suggestions.

\end{document}